\documentclass[letterpaper, 10 pt, conference]{ieeeconf}  

\IEEEoverridecommandlockouts                              

\usepackage{amsmath} 
\usepackage{amssymb}  
\usepackage{booktabs}

\usepackage{graphicx}
\title{\LARGE \bf
Introducing the quadratically-constrained quadratic programming framework in HPIPM
}

\author{Gianluca Frison$^{1}$, Jonathan Frey$^{1}$, Florian Messerer$^{1}$, Andrea Zanelli$^{2}$ and Moritz Diehl$^{1}$
\thanks{*This research was supported by DFG via Research Unit FOR 2401 and project 424107692 and by the EU via ELO-X 953348 and by INTERREG V Upper Rhine, project ACA-MODES.}
\thanks{$^{1}$Department of Microsystems Engineering (IMTEK), University of Freiburg
        {\tt\small \{gianluca.frison, jonathan.frey, florian.messerer, moritz.diehl\} @imtek.uni-freiburg.de}}%
\thanks{$^{2}$ Institute for Dynamic Systems and Control (IDSC), ETH Zurich
        {\tt\small zanellia@ethz.ch}}%
}

\begin{document}

\makeatletter
\renewcommand*\env@matrix[1][*\c@MaxMatrixCols c]{%
  \hskip -\arraycolsep
  \let\@ifnextchar\new@ifnextchar
  \array{#1}}
\makeatother

\maketitle
\thispagestyle{empty}
\pagestyle{empty}

\begin{abstract}

This paper introduces the quadratically-constrained quadratic programming (QCQP) framework recently added in HPIPM alongside the original quadratic-programming (QP) framework.
The aim of the new framework is unchanged, namely providing the building blocks to efficiently and reliably solve (more general classes of) optimal control problems (OCP).
The newly introduced QCQP framework provides full features parity with the original QP framework: three types of QCQPs (dense, optimal control and tree-structured optimal control QCQPs) and interior point method (IPM) solvers as well as (partial) condensing and other pre-processing routines.
Leveraging the modular structure of HPIPM, the new QCQP framework builds on the QP building blocks and similarly provides fast and reliable IPM solvers.

\end{abstract}

\section{Introduction}

In the field of optimal control, Quadratic Programs (QP) are well established and arise in many contexts, both directly, e.g., in linear Model Predictive Control (MPC), and indirectly as subproblems in Sequential Quadratic Programming (SQP) algorithms for Nonlinear Programming (NLP).
In consequence, many efficient solver implementations exist.
However, as QPs admit only affine constraints, they lack some important expressiveness.
Quadratic inequality constraints occur naturally in many problems, e.g., in MPC for power electronics \cite{Almer2021}, electric drives \cite{Geweth2020, Zanelli2021a}, or ellipsoidal corridors in robotics \cite{Duijkeren2019}, to name just a few.
Furthermore, MPC formulations with desiderable theoretical properties such as stability guarantees via ellipsoidal terminal sets~\cite{Mayne2000a} or robustness against modeling errors and noise in a tube MPC framework~\cite{Mayne2005} define quadratic inequality constraints alike.
If the equality constraints remain linear, this leads to the class of Quadratically Constrained Quadratic Programs (QCQP) \cite{Boyd2004}.
For some of the mentioned applications, the sampling times are in the three digit microsecond range, creating the need for highly efficient solvers.
One approach leading to satisfactory results in praxis may be to approximate the ellipsoidal constraints by an inner polytope via affine inequalities in order to obtain a QP. As we will demonstrate on a numerical example, already for a moderate number of affine constraints it can be faster to directly solve the QCQP.

Since roughly the 2000s, some interest has developed in using Sequential QCQP (SQCQP) methods for NLPs, which sequentially solve QCQP approximations, see \cite{Izmailov2014} for an overview.
While this idea is certainly older, it was initially not of practical interest due to the lack of efficient solvers for the QCQP subproblems \cite{Fukushima1986}.
Intuitively, QCQP can provide better approximations of an NLP than a QP would, as they are able to capture the constraint curvature also beyond their Hessian contribution at the linearization point.
Therefore, one would expect SQCQP methods to require fewer iterations than SQP methods, although each iteration will be more expensive.
However, if an efficient QCQP solver is employed, this may tip the balance in favor of SQCQP.
Additionally, SQCQP methods can partially avoid the Maratos effect \cite{Maratos1978}, which is a consequence of the constraint linearization inherent to SQP methods \cite{Fukushima2003a}.
If the NLP is generally nonconvex but exhibits convex-over-nonlinear structures, convexity of the QCQP subproblems may be assured via generalized Gauss-Newton Hessian approximations \cite{Messerer2021a}.
For quadratic convexities, this is equivalent to Sequential Convex Programming (SCP), cf. \cite{Messerer2020, TranDinh2010}, which in this case also leads to convex QCQP subproblems.
While QCQPs are often formulated and solved as Second Order Cone Programs (SOCP) \cite{Boyd2004}, we believe and illustrate via numerical experiments that important efficiency gains can be made by treating them directly as QCQP.

Besides the availability of general large scale sparse SOCP solvers, the interest over QCQPs has led to the development of some software packages specifically targeting embedded applications, and possibly tailored to the optimal control problem structure.
ECOS~\cite{Domahidi2013a} is a sparse SOCP solver targeting embedded applications.
The commercial software FORCES~\cite{Domahidi2012} provides a QCQP solver tailored to multistage problems arising in optimal control.
More recently, in \cite{Chen2020} two algorithms for QCQP are presented, which are built on the basis of the open-source QP solvers OSQP \cite{Stellato2020} and HPIPM \cite{Frison2020a} respectively, although tailored to a specific MPC problem.

In this paper, we present an extension of the convex QP framework of HPIPM \cite{Frison2020a} to the class of convex QCQP.
This choice, in contrast to the more general class of SOCP, allows us to provide full feature parity with the existing QP framework, e.g. supporting (partial) condensing for QCQP too.
The HPIPM software package has been extended by building blocks for the efficient solution of QCQP-based problems in MPC, similar to the existing QP framework.
The QCQP formulation and framework on purpose resemble their QP version to provide a familiar interface to existing users and higher-level algorithms.

\section{QCQP formulations}

As a direct extension to the QP framework described in~\cite{Frison2020a}, HPIPM currently defines three QCQP types: a dense QCQP, an OCP QCQP and tree-structured OCP QCQP.
As a distinctive feature inherited from the QP framework, all QCQP formulations define a special type of variable, the slacks, which do not enter the dynamics equality constraints and which give a diagonal contribution to Hessian and inequality constraint matrices.
The exploitation of this structure allows their elimination from the QCQP formulation in computational complexity linear in the number of slack variables, 
making it computationally cheap to use them.
The slack variables can be used to efficiently implement soft constraints with L1 and L2 penalties on all constraint types, comprising the quadratic constraints, as done extensively in Section~\ref{sec:num_res}.

\subsection{Dense QCQP}

The dense QCQP type describes a generic QCQP where Hessian and constraint matrices are assumed to be dense. 
This formulation can handle the QCQP sub-problems arising in single-shooting discretization schemes or in state-condensing schemes in OCP and MPC frameworks.
The formulation is stated as
\begin{align*}
\min_{v,s} & \quad \frac 1 2 \begin{bmatrix} v \\ 1 \end{bmatrix}^T \begin{bmatrix} H & g \\ g^T & 0 \end{bmatrix} \begin{bmatrix} v \\ 1 \end{bmatrix} \\
& + \frac 1 2 \begin{bmatrix} s^{\rm l} \\ s^{\rm u} \\ 1 \end{bmatrix}^T \begin{bmatrix} Z^{\rm l} & 0 & z^{\rm l} \\ 0 & Z^{\rm u} & z^{\rm u} \\ (z^{\rm l})^T & (z^{\rm u})^T & 0 \end{bmatrix} \begin{bmatrix} s^{\rm l} \\ s^{\rm u} \\ 1 \end{bmatrix} \\
{\rm s.t.} & \quad A v = b \\
& \quad \begin{bmatrix} \underline v \\ \underline d \end{bmatrix} \leq \begin{bmatrix} J^{b,v} \\ C \end{bmatrix} v + \begin{bmatrix} J^{s,v} \\ J^{s,g} \end{bmatrix} s^{\rm l} \\
& \quad \begin{bmatrix} J^{b,v} \\ C \end{bmatrix} v - \begin{bmatrix} J^{s,v} \\ J^{s,g} \end{bmatrix} s^{\rm u} \leq \begin{bmatrix} \overline v \\ \overline d \end{bmatrix} \\
& \quad \frac 1 2 \begin {bmatrix} v \\ 1 \end{bmatrix}^T \begin{bmatrix} H_k & g_k \\ g_k^T & 0 \end{bmatrix} \begin{bmatrix} v \\ 1 \end{bmatrix} - J^{s,q}_k s^{\rm u} \leq \overline d_k \; , \; k \in \mathcal {H}_q \\
& \quad s^{\rm l}\geq \underline s^{\rm l} \\
& \quad s^{\rm u}\geq \underline s^{\rm u}
\end{align*}
where the quadratic constraint index $k$ takes values in the set $\mathcal {H}_q=\{1,\dots,n_\mathrm{q}\}$.
The dense QCQP formulation closely resembles the dense QP formulation with the addition of quadratic constraints.
The primal variables comprise generic variables $v$ and slack variables $s^{\rm l}$ ($s^{\rm u}$) associated to the lower (upper) constraints.
The Hessian matrices of the slacks $Z^{\rm l}$ and $Z^{\rm u}$ are diagonal.
The matrices $J^{\cdot , \cdot}$ are made of rows from identity matrices, and are employed to select only some components in box and soft constraints.
The constraint matrices are the same for the upper and the lower constraints, meaning that all linear constraints in the formulation are two-sided.
On the other hand, the quadratic constraints are only one-sided in order to ensure convexity, since the matrices $H_k$ are assumed to be positive semi-definite.
Also, they are scalar, i.e. $\bar d_k$ is a real number and $J_k^{s,q}$ a row vector.
A mask (not represented in the above dense QCQP formulation) can be employed to dynamically activate or deactivate the single upper and/or lower constraints.
The considerations in this paragraph apply also to the OCP and tree OCP QCQP.

\subsection{Optimal Control Problem (OCP) QCQP}

The OCP QCQP type describes a QCQP formulation handling many OCP and MPC problems such as constrained linear MPC problems and QCQP sub-problems in SQCQP and SCP algorithms for non-linear OCP and MPC problems~\cite{Messerer2021a}.
The formulation is stated as
\begin{align*}
\min_{x,u,s} & \quad \sum_{n=0}^N \frac 1 2 \begin{bmatrix} u_n \\ x_n \\ 1 \end{bmatrix}^T \begin{bmatrix} R_n & S_n & r_n \\ S_n^T & Q_n & q_n \\ r_n^T & q_n^T & 0 \end{bmatrix} \begin{bmatrix} u_n \\ x_n \\ 1 \end{bmatrix} \\
& + \frac 1 2 \begin{bmatrix} s^{\rm l}_n \\ s^{\rm u}_n \\ 1 \end{bmatrix}^T \begin{bmatrix} Z^{\rm l}_n & 0 & z^{\rm l}_n \\ 0 & Z^{\rm u}_n & z^{\rm u}_n \\ (z^{\rm l}_n)^T & (z^{\rm u}_n)^T & 0 \end{bmatrix} \begin{bmatrix} s^{\rm l}_n \\ s^{\rm u}_n \\ 1 \end{bmatrix} \\
{\rm s.t.} & \quad x_{n+1} = A_n x_n + B_n u_n + b_n \quad , \quad n \in \mathcal H \backslash \{N\} \\ 
& \quad \begin{bmatrix} \underline u_n \\ \underline x_n \\ \underline d_n \end{bmatrix} \leq \begin{bmatrix} J_n^{b,u} & 0 \\ 0 & J_n^{b,x} \\ D_n & C_n \end{bmatrix} \begin{bmatrix} u_n \\ x_n \end{bmatrix} + \begin{bmatrix} J_n^{s,u} \\ J_n^{s,x} \\ J_n^{s,g} \end{bmatrix} s^{\rm l}_n \; , \; n \in \mathcal H \\ 
& \quad \begin{bmatrix} J_n^{b,u} & 0 \\ 0 & J_n^{b,x} \\ D_n & C_n \end{bmatrix} \begin{bmatrix} u_n \\ x_n \end{bmatrix} - \begin{bmatrix} J_n^{s,u} \\ J_n^{s,x} \\ J_n^{s,g} \end{bmatrix} s^{\rm u}_n \leq \begin{bmatrix} \overline u_n \\ \overline x_n \\ \overline d_n \end{bmatrix} \; , \; n \in \mathcal H \\ 
& \quad \frac 1 2 \begin{bmatrix} u_n \\ x_n \\ 1 \end{bmatrix}^T \begin{bmatrix} R_{n,k} & S_{n,k} & r_{n,k} \\ S_{n,k}^T & Q_{n,k} & q_{n,k} \\ r_{n,k}^T & q_{n,k}^T & 0 \end{bmatrix} \begin{bmatrix} u_n \\ x_n \\ 1 \end{bmatrix} \\
& \qquad \qquad \qquad \qquad \quad - J_{n,k}^{s,q} s_n^{\rm u} \leq \overline u_{n,k} \quad , \quad \left\{ \begin{matrix} n \in \mathcal H \\ k \in \mathcal H_q \end{matrix} \right. \\
& \quad s_n^{\rm l}\geq \underline s_n^{\rm l} \; , \; n \in \mathcal H \\ 
& \quad s_n^{\rm u}\geq \underline s_n^{\rm u} \; , \; n \in \mathcal H 
\end{align*}
where the stage index $n$ takes values in the set $\mathcal H = \{0,1,\dots,N\}$ and the quadratic constraint index $k$ takes values in the set $\mathcal {H}_q=\{1,\dots,n_\mathrm{q}\}$.
This problem has a multi-stage structure, with cost and inequality constraints defined stage-wise, and with dynamics equality constraints coupling pairs of consecutive stages.
All data matrices and vectors can vary stage-wise.
The primal variables are divided into state variables $x_n$, control (or input) variables $u_n$ and slack variables associated to the lower (upper) constraints $s^{\rm l}_n$ ($s^{\rm u}_n$).
The size of all variables (number of states $n_{x_n}$, number of controls $n_{u_n}$ and number of slacks $n_{s_n}$), as well as the number of box constraints $n_{b_n}$, general polytopic constraints $n_{g_n}$ and quadratic constraints $n_{q_n}$ can freely vary stage-wise, as required e.g. to solve QCQP sub-problems arising in multi-phase multiple-shooting discretization schemes, similarly to QP sub-problems in~\cite{Bock1984}.

Note that the current formulation does not explicitly define equality constraints other than the dynamics equations, and therefore other types of equality constraints have to be reformulated as inequality constraints with equal upper and lower limits.
Currently the HPIPM framework provides some functionality (not represented in the above OCP QCQP formulation) to explicitly mark which of the inequality constraints should be considered as equalities, and routines to remove simple types of equality constraints (as e.g. the constraint on the initial state value $x_0 = \hat x_0$) from the problem formulation in a pre-solve step.

\subsection{Tree OCP QP}

The tree OCP QCQP type can handle a large number of common robust and scenario-based OCP and MPC problems, see~\cite{Kouzoupis2019} and references therein.
The formulation is stated as
\begin{align*}
\min_{x,u,s} & \quad \sum_{n \in \mathcal N} \frac 1 2 \begin{bmatrix} u_n \\ x_n \\ 1 \end{bmatrix}^T \begin{bmatrix} R_n & S_n & r_n \\ S_n^T & Q_n & q_n \\ r_n^T & q_n^T & 0 \end{bmatrix} \begin{bmatrix} u_n \\ x_n \\ 1 \end{bmatrix} \\
& + \frac 1 2 \begin{bmatrix} s^{\rm l}_n \\ s^{\rm u}_n \\ 1 \end{bmatrix}^T \begin{bmatrix} Z^{\rm l}_n & 0 & z^{\rm l}_n \\ 0 & Z^{\rm u}_n & z^{\rm u}_n \\ (z^{\rm l}_n)^T & (z^{\rm u}_n)^T & 0 \end{bmatrix} \begin{bmatrix} s^{\rm l}_n \\ s^{\rm u}_n \\ 1 \end{bmatrix} \\
{\rm s.t.} & \quad x_m = A_m x_n + B_m u_n + b_m \quad , \quad \left\{ \begin{matrix} n \in \mathcal N \backslash \mathcal L \\ m \in \mathcal C(n) \end{matrix} \right. \\
& \quad \begin{bmatrix} \underline u_n \\ \underline x_n \\ \underline d_n \end{bmatrix} \leq \begin{bmatrix} J_n^{b,u} & 0 \\ 0 & J_n^{b,x} \\ D_n & C_n \end{bmatrix} \begin{bmatrix} u_n \\ x_n \end{bmatrix} + \begin{bmatrix} J_n^{s,u} \\ J_n^{s,x} \\ J_n^{s,g} \end{bmatrix} s^{\rm l}_n \; , \; n \in \mathcal N \\
& \quad \begin{bmatrix} J_n^{b,u} & 0 \\ 0 & J_n^{b,x} \\ D_n & C_n \end{bmatrix} \begin{bmatrix} u_n \\ x_n \end{bmatrix} - \begin{bmatrix} J_n^{s,u} \\ J_n^{s,x} \\ J_n^{s,g} \end{bmatrix} s^{\rm u}_n \leq \begin{bmatrix} \overline u_n \\ \overline x_n \\ \overline d_n \end{bmatrix} \; , \; n \in \mathcal N \\
& \quad \frac 1 2 \begin{bmatrix} u_n \\ x_n \\ 1 \end{bmatrix}^T \begin{bmatrix} R_{n,k} & S_{n,k} & r_{n,k} \\ S_{n,k}^T & Q_{n,k} & q_{n,k} \\ r_{n,k}^T & q_{n,k}^T & 0 \end{bmatrix} \begin{bmatrix} u_n \\ x_n \\ 1 \end{bmatrix}\\
& \qquad \qquad \qquad \qquad \quad - J_{n,k}^{s,q} s_n^{\rm u} \leq \overline u_{n,k} \quad , \quad \left\{ \begin{matrix} n \in \mathcal N \\ k \in \mathcal H_q \end{matrix} \right. \\
& \quad s_n^{\rm l}\geq \underline s_n^{\rm l} \; , \; n \in \mathcal N \\
& \quad s_n^{\rm u}\geq \underline s_n^{\rm u} \; , \; n \in \mathcal N
\end{align*}
where $\mathcal N$ is the set of nodes in the tree and $\hat N$ is its cardinality.
The set $\mathcal L$ contains the leaves of the tree, while $\mathcal C(n)$ denotes the set of the children of node $n$.
All data matrices and vector as well as their size can vary node-wise.

\section{Algorithm}

As outlined in~\cite{Frison2020a}, HPIPM provides all algorithmic building blocks needed to implement many variants of a primal-dual IPM in a modular fashion.
In the current paper, these building blocks are used to implement efficient algorithms for QCQP.

\subsection{Primal-dual IPM}

This section contains some basic notion about primal-dual IPM algorithms, with the only purpose of giving some intuition on the algorithm derivation.
It is implicitly assumed that all the necessary regularity conditions and constraint qualifications hold.
The algorithm derivation is at first targeting a generic optimization problem, and is subsequently tailored to the QCQP case.
The IPM algorithm is derived from the KKT conditions, and not from the barrier methods theory.
The interested reader can find more detailed and rigorous presentations in~\cite{nocedal2006, wright1997}.

Let us consider an optimization problem (OP) in the form
\begin{equation}
\begin{aligned}
\min_y & \quad f(y) \\
s.t. & \quad g(y) = 0 \\
& \quad h(y) \geq 0.
\end{aligned}
\label{eq:alg:pdipm:op}
\end{equation}
The Lagrangian function for this OP reads
\begin{equation*}
\mathcal L (y, \pi, \lambda) = f(y) - \pi^T g(y) - \lambda^T h(y)
\end{equation*}
where $\pi$ and $\lambda$ are the Lagrange multipliers of the equality and inequality constraints, respectively.
The first order necessary KKT optimality conditions read
\begin{subequations}
\begin{align}
& \nabla_y \mathcal L(y, \pi, \lambda) = \nabla f(y) - \nabla g(y) \pi - \nabla h(y) \lambda = 0 \label{eq:alg:pdipm:kkt:g} \\
& - g(y) = 0 \label{eq:alg:pdipm:kkt:b} \\
& - h(y) + t = 0 \label{eq:alg:pdipm:kkt:d} \\
& \lambda_i t_i = 0 \qquad \qquad \qquad \qquad \qquad \qquad i=1,\dots,n_\mathrm{i} \label{eq:alg:pdipm:kkt:m} \\
& (\lambda, t) \geq 0 \label{eq:alg:pdipm:kkt:s}
\end{align}
\label{eq:alg:pdipm:kkt}
\end{subequations}
where the slack variables $t = h(y) \geq 0$ have been introduced, and $n_\mathrm{i}$ denotes the total number of inequality constraints.

Equations (\ref{eq:alg:pdipm:kkt:g})-(\ref{eq:alg:pdipm:kkt:m}) are a system of nonlinear equations $F(y, \pi, \lambda, t)=0$. 
In a nutshell, a primal-dual IPM is Newton method applied to the system of equations $F_{\mu}(y, \pi, \lambda, t)=0$ with equation (\ref{eq:alg:pdipm:kkt:m}) relaxed as
\begin{equation*}
\lambda_i t_i = \mu, \quad i=1,\dots,n_\mathrm{i}.
\end{equation*}
The homotopy parameter $\mu$ is related to the barrier parameter in barrier methods, and it is shrunk toward zero as the iterations approach the solution of (\ref{eq:alg:pdipm:kkt}).
A line search procedure is used to ensure the strict satisfaction of the inequalities (\ref{eq:alg:pdipm:kkt:s}) on the sign of Lagrange multipliers and slacks of inequality constraints.
Furthermore, in general a line search or a trust-region algorithm can be employed to ensure a sufficient progress along the search direction (although this is not employed in the HPIPM implementation described in the current paper).

At every iteration $k$ of the Newton method, the Newton step $(\Delta y_{\rm aff}, \Delta \pi_{\rm aff}, \Delta \lambda_{\rm aff}, \Delta t_{\rm aff})$ is found by solving the linear system
\begin{multline*}
\nabla F_{\mu}(y_k, \pi_k, \lambda_k, t_k) \begin{bmatrix} \Delta y_{\rm aff}^T & \Delta \pi_{\rm aff}^T & \Delta \lambda_{\rm aff}^T & t_{\rm aff}^T \end{bmatrix}^T \\
= - F_{\mu}(y_k, \pi_k, \lambda_k, t_k)
\end{multline*}
which takes the form
\begin{multline}
\!\!\!\!\!\!\!\!\begin{bmatrix} \nabla_{yy}^2 \mathcal L (y_k,\pi_k,\lambda_k) & - \nabla g(y_k) & - \nabla h(y_k) & 0 \\ - (\nabla g(u_k))^T & 0 & 0 & 0 \\ - (\nabla h(u_k))^T & 0 & 0 & I \\ 0 & 0 & T_k & \Lambda_k \end{bmatrix}
\!\!\begin{bmatrix} \Delta y_{\rm aff} \\ \Delta \pi_{\rm aff} \\ \Delta \lambda_{\rm aff} \\ \Delta t_{\rm aff} \end{bmatrix} \\
= - \begin{bmatrix} \nabla f(y_k) - \nabla g(y_k) \pi_k - \nabla h(y_k) \lambda_k \\ - g(y_k) \\ - h(y_k)+ t_k \\ \Lambda_k T_k e - \mu e \end{bmatrix}  =: - \begin{bmatrix} r_g \\ r_b \\ r_d \\ r_m \end{bmatrix}
\label{eq:alg:pdipm:op_delta}
\end{multline}
where
\begin{align}
&\nabla_{yy}^2 \mathcal L (y_k,\pi_k,\lambda_k) \\
&= \nabla^2 f(y_k) - \sum_{i=0}^{n_\mathrm{e}-1} \pi_{k,i} \nabla^2 g_i(y_k) - \sum_{i=0}^{n_{\mathrm i}-1} \lambda_{k,i} \nabla^2 h_i(y_k) \nonumber
\end{align}
and where $n_\mathrm{e}$ is the number of equality constraints.
The diagonal matrices $\Lambda_k$ and $T_k$ have on their diagonal the elements of the vector $\lambda_k$ and $t_k$ respectively.
The function $F_{\mu}(y_k, \pi_k, \lambda_k, t_k)$ at the right hand side (RHS) is denoted as the residual function.

In order to tailor the results to the QCQP case, let us consider the simple QCQP with linear equality constraints and only quadratic inequality constraints (the linear inequality constraints being a special case)
\begin{align}
\min_y & \quad \frac 1 2 y^T \mathcal H_0 y + g_0^T y \nonumber \\
\mathrm{s.t.} & \quad \mathcal A y - b = 0 \label{eq:alg:pdipm:qcqp} \\
& \quad - \frac 1 2 y^T \mathcal H_i y - g_i^T y + d_i \geq 0 \quad , \quad i=\{1,\dots,n_\mathrm{q}\} \nonumber
\end{align}

The Newton step is found solving the linear system 
\begin{multline}
\begin{bmatrix} \mathcal H(\lambda_k) & - \mathcal A^T & - \mathcal C(y_k)^T & 0 \\ - \mathcal A & 0 & 0 & 0 \\ - \mathcal C(y_k) & 0 & 0 & I \\ 0 & 0 & T_k & \Lambda_k \end{bmatrix}
\begin{bmatrix} \Delta y_{\rm aff} \\ \Delta \pi_{\rm aff} \\ \Delta \lambda_{\rm aff} \\ \Delta t_{\rm aff} \end{bmatrix}  \\
= - \begin{bmatrix} \mathcal H_0 y_k - \mathcal A^T \pi_k - \mathcal C(y_k)^T \lambda_k + g_0 \\ - \mathcal A y_k + b\\ - \mathcal C(y_k) y_k + t_k + d \\ \Lambda_k T_k e - \mu e \end{bmatrix} =: - \begin{bmatrix} r_g \\ r_b \\ r_d \\ r_m \end{bmatrix}
\label{eq:alg:pdipm:qcqp_delta}
\end{multline}
where
\begin{equation*}
\begin{aligned}
\mathcal H(\lambda_k) &= \mathcal H_0 + \sum_{i=1}^{n_\mathrm{q}} \lambda_{k,i} \mathcal H_i \\
\mathcal C(y_k) &= \begin{bmatrix} - \mathcal H_1 y_k - g_1 \\ \vdots \\ - \mathcal H_{n_\mathrm{q}} y_k - g_{n_\mathrm{q}} \end{bmatrix} \quad , \quad 
d = \begin{bmatrix} - d_1 \\ \vdots \\ - d_{n_\mathrm{q}} \end{bmatrix}
\end{aligned}
\end{equation*}
The linear system~(\ref{eq:alg:pdipm:qcqp_delta}) resembles the one solved in the QP case in~\cite{Frison2020a}, with the notable differences that the matrices $\mathcal H(\lambda_k)$ and $\mathcal C(y_k)$ depend on the value of the variables $y_k$ and $\lambda_k$ at iteration $k$.
Therefore, an IPM for a QCQP can be implemented by using the building block provided by the same linear system factorization already present in the QP framework in HPIPM, provided that the matrices $\mathcal H(\lambda_k)$ and $\mathcal C(y_k)$ are updated beforehand at each iteration $k$.

\subsection{Implementation choices}

This section discusses some of the main choices in implementing the QCQP framework in HPIPM.
Most of these choices are directly inherited from the QP framework so are only briefly mentioned here and mostly to highlight possible differences or similarities.
More details can be found in the original QP framework description in~\cite{Frison2020a}.

\subsubsection{Barrier parameter and step length selection}

Similarly to the QP IPM solver, the barrier parameter $\mu$ is chosen adaptively using a Mehrotra's predictor-corrector based strategy.
More precisely, this is implemented using the so called conditional predictor-corrector variant proposed in~\cite{nocedal2009}, in which the corrector step is employed only if it does not increase the duality measure too much, therefore safeguarding against its possible harmful effects.
The conditional Mehrotra's predictor-corrector is shown to provide robust convergence also in case of nonlinear programming and trivial initial point strategies, even when globalization strategies are not employed~\cite{nocedal2009}.

The primal and dual step lengths are computed as the longest steps such that the constraint on the sign of the slacks and the Lagrange multipliers of the inequality constraints are not violated.
Two different strategies, one based on separate steps for the primal and the dual variables and one based on a single step for all variables can be employed depending on the speed/robustness trade-off.
An adaptive fraction to the boundary parameter $\tau$ is employed, which is chosen in a conservative fashion when the step is short, while it allows for a fast progress towards the solution when the step is close to a full step.
Thanks to the robustness of the conditional Mehrotra's predictor-corrector strategy, further globalization strategies are currently not implemented, but they will be investigated in future research.


\subsubsection{Newton system solution and iterative refinement}

The Newton system solution and the optional iterative refinement steps are directly inherited from the QP framework and leverage the efficient building blocks implemented therein.

In summary, the slacks $\Delta t_{\rm aff}$ and the Lagrange multipliers of the inequality constraints $\Delta \lambda_{\rm aff}$ are eliminated before the actual factorization is performed: this takes the linear system in the form of the KKT system of an equality constrained QP (usually named `augmented system'), allowing efficient factorization procedures where the pivot sequence is fixed and the problem structure is fully exploited using computations on the dense sub-matrices.
The drawback of this procedure is that the elimination of slacks and Lagrange multipliers causes ill-conditioning, since as the iterates approach the solution some of their components converge to finite quantities while others converge to zero, and therefore their ratio diverges.

All currently implemented QP types are treated using dense or structured (where the overall sparsity is exploited by a hand-crafted algorithm operating on dense sub-matrices) factorization procedures. No arbitrary sparsity patterns whithin the problem data are exploited.
In the dense QCQP case, the augmented system is solved using either null-space or Schur-complement methods~\cite{nocedal2006}. In the case of fully condensed MPC problems there are generally no equality constraints, so simply a Cholesky factorization is used.
In the case of the OCP QCQP, the augmented system is factorized using either a classic or a square-root Riccati recursion~\cite{frison2015a}, possibly implemented using QR-based array algorithms for improved numerical accuracy and stability; the factorization has a computational complexity of $\mathcal O(N(n_x+n_u)^3)$ flops.
In the case of the tree OCP QCQP, the augmented system is factorized using a Riccati recursion modified to exploit the tree structure without introducing fill-in outside the data matrices~\cite{frison2017a}. All Riccati variants in the OCP QCQP type have a corresponding tree-tailored variant, and their computational complexity is of $\mathcal O(\hat N(n_x+n_u)^3)$ flops.

In case of unstable or ill-conditioned systems, especially at late IPM iterations, it can happen that the accuracy of the Newton step is too low to be useful for the IPM algorithm.
In such cases, it may be useful to perform a few iterative refinement steps, where the same linear system factorization is employed to solve for the residuals of the linear system to compute a correction term and iteratively improve the accuracy of the Newton step.
Note however that in case the numerical factorization is a too poor approximation of the exact one, iterative refinement may get unstable and further degrade the accuracy of the solution, so it should be used with care.
Iterative refinement can also be useful to compensate for the effect of regularization, that may be needed to successfully complete the linear system factorization.

\subsubsection{Delta and absolute IPM formulations}

In the QP case, there is a computational advantage in reformulating the linear system~\eqref{eq:alg:pdipm:op_delta} such that it directly computes the new iterate instead of the step.
This exploits the linearity of equations (\ref{eq:alg:pdipm:kkt:g})-(\ref{eq:alg:pdipm:kkt:d}) in the QP case by replacing the residuals computation at the RHS (which has a quadratic cost in the number of stage variables) with a vector difference (which has a linear cost in the number of stage variables), at the cost of increased numerical instability and possibly severe cancelation errors at late IPM iterations.
Nonetheless, for small and well conditioned QPs in applications where only low accuracy is required, this can give a noticeable performance boost.

In the QCQP case, equations \eqref{eq:alg:pdipm:kkt:g} and \eqref{eq:alg:pdipm:kkt:d} are not linear any longer, so the computational advantage of using an absolute IPM formulation is strongly reduced, while the drawbacks remain.
Therefore, the absolute formulation is not considered any longer in this paper and all algorithmic variants are based on the standard delta formulation.

\subsubsection{Other implementation choices}

The standard infeasible-start Mehrotra's predictor-corrector IPM algorithm is tweaked by several additional parameters whose different choice gives rise to a family of algorithms with different trade-offs between computational speed and reliability.
In particular, a `mode' argument is introduced in order to select pre-defined sets of parameters tailored to a particular purpose.
In the QCQP solver, the {\tt speed}, {\tt balance} and {\tt robust} modes are available, while the {\tt speed\_abs} is not, as previously discussed. All modes work similarly to their analogues in the QP solver.
Just to mention one notable feature, in the {\tt balance} mode, in case the accuracy of the linearized KKT system factorization is too low, the factorization is repeated by replacing all Cholesky factorizations of normal matrices in the form $A\cdot A^T$ by array algorithms based on QR factorizations, which have better numerical properties as they never explicitly form the worse-conditioned normal matrix $A\cdot A^T$.
In the {\tt robust} mode, the more accurate (but slower) QR factorization based array algorithms are always employed.
See~\cite{Frison2020a} for a more descriptive list of other implementation choices directly inherited from the QP framework.

\subsubsection{Dense linear algebra}

Algorithms in HPIPM are explicitly designed to exploit the overall block-sparse structure of the optimization problems.
Sparsity within data and working matrices is not considered.
Therefore, computations are cast in terms of operations on matrices with fixed structure, like e.g. dense, symmetric, triangular or diagonal.

In HPIPM, the KKT and the (partial) condensing modules contain the most computationally expensive routines, which have cubic (linearized KKT system factorization, condensing of LHS) or quadratic (linearized KKT system solution, residuals computation, condensing of RHS) complexity in the number of (stage) variables.
In case of all currently implemented QP and QCQP types, these modules are implemented using the structure-based interface in the high-performance linear algebra library BLASFEO~\cite{Frison2018}. 
These routines are designed to give close-to-peak performance for matrices of moderate size (assumed to fit in cache) and are optimized for several hardware architectures (e.g. exploiting different classes of SIMD instructions).

All other HPIPM modules are completely self-contained and independent of the matrix type or linear algebra library employed in the KKT modules.

\subsection{(Partial) condensing of quadratic constraints}

One key advantage of directly targeting QCQPs instead of reformulating them as more general problems such as SOCPs is the possibility to extend in a straightforward way all the rich and efficient existing framework for QPs.
In particular, one powerful tool in the HPIPM QP framework is certainly state condensing, both in the full and the partial condensing variants.
In a nutshell, full condensing exploits the dynamics equality constraints to remove all superflous state variables and only keep the true degrees of freedom in the condensed problem formulation, namely the initial state (if not previously eliminated) and the control variables.
Partial condensing is a more recently developed technique~\cite{Axehill2015} that performs condensing in $N_c$ blocks of consecutive stages, resulting in a problem that still has the form of an OCP but with shorter horizon $N_c$ and larger control vector size per stage.
The size of the blocks (and therefore the horizon length of the partially condensed problem $N_c$) is a tuning parameter that can be used to trade off horizon length with control vector size to e.g. minimize the flop count or the solution time for every specific QP solver.
Furthermore, the size of the dense sub-matrices in the OCP formulation is increased, generally resulting in improved performance of the dense basic linear algebra routines~\cite{frison2016}.
In case of OCPs with many more states than controls and rather long horizon (which is a fairly common case in practice), partial condensing can result in multiple times faster solution times compared to both the original full space formulation and the fully condensed one.

In the QCQP case, the only new bit is the condensing of the quadratic constraints.
This is analogue to the condensing of the quadratic cost, with the key difference that the cost is a sum of quadratic terms (one per stage), while each quadratic constraint is a single quadratic function of controls and states at a certain stage.
Therefore, out of the three possible Hessian condensing algorithms~\cite{frison2015a}, the `classical' one (i.e. the one with complexity $\mathcal O(N^3)$ and $\mathcal O(n_x^2)$) is the most efficient variant in case of quadratic constraints, since it allows one to only perform the computations for that certain stage and it does not require recursions involving other stages.
The computational cost to condense a quadratic constraint depends on the specific stage $n$ the constraint belongs to.
It is of
\begin{equation*}
n^2 n_u^2 n_x + 2 n n_u n_x^2
\end{equation*}
flops if the initial state is not an optimization variable, and
\begin{equation*}
n^2 n_u^2 n_x + 4 n n_u n_x^2 + 3 n_x^3
\end{equation*}
flops if the initial state is an optimization variable (such as in partial condensing of middle blocks).
Notice that the cost increases with the stage index $n$ and it is therefore largest at the last stage where $n=N$.
As a comparison, the cumulative cost of condensing one quadratic constraint per stage is roughly equivalent to the cost of condensing the Hessian using the `classical' algorithm.
This cost can be reduced in the (rather common) case of diagonal or zero quadratic constraint matrices, therefore some functionality to handle such cases efficiently is implemented in HPIPM.

\section{Numerical results} \label{sec:num_res}

The dynamical system used throughout in this numerical results section is the linear mass-spring system described in~\cite{wang2010a}.
More precisely, we consider a linear chain of masses of 1 $ \mathrm{kg} $ connected to each other and to the walls on either side with springs of spring constant $ 1 \mathrm{N}{\rm m}^{-1} $.
The states are the relative position and the velocity of each mass, and therefore $n_x$ is twice the number of masses.
As controls we consider the actuator forces acting on the first $n_u$ masses.
This simple example allows one to easily scale the number of states, controls and the horizon length, as well as imposing different types of constraints: for example it is straightforward to define quadratic constraints that describe the total energy (sum of masses kinetic and springs potential energies) of the system.

All numerical experiments are run on a laptop equipped with a Intel Core i7 4810MQ processor.


\subsection{Comparison of QP and QCQP frameworks in HPIPM}

%

In this section, we compare the QP and the QCQP frameworks in HPIPM in order to analyze the intrisic efficiency differences between equivalently implemented QP and QCQP solvers.

In Table~\ref{tab:qp_qcqp:1}, there is the comparison of a QP with two QCQP where a different number of box constraints is replaced with quadratic constraints.
All problems have $N=15$ horizon, $n_x=4$ states, $n_u=1$ controls; additionally they have:
\begin{itemize}
\item QP\_0: 1 control bound per stage; 4 terminal softed state bounds
\item QCQP\_1: same 1 control bound per stage; 1 terminal softed state quadratic constraint (replacing the 4 terminal softed state bounds)
\item QCQP\_N: 1 control quadratic constraint per stage (replacing the 1 control bound per stage); same 1 terminal softed state quadratic constraint
\end{itemize}
The solution times for the QCQPs are within roughly a factor 2 of the solution time for the QP for all algorithmic combinations.
QCQP with (partial) condensing shows some additional overhead compared to the QP counterparts but still gives very noticeable speedups, in excess of a factor 2 for these problems.
In QCQP problems with only one terminal quadratic constraint, the overhead compared to an analogue QP is rather limited and in the range of 15-30\% for the considered problems.

\begin{table}
\centering
\caption{Total solution times (in seconds, for fixed 7 IPM iterations plus $x_0$ removal and condensing) for the mass-spring QP and QCQPs with different number of quadratic constraints}
\label{tab:qp_qcqp:1}
\begin{tabular}{@{}llll@{}}
\toprule
problem & QP\_0 & QCQP\_1 & QCQP\_N \\
\midrule
baseline                  & $9.82\cdot 10^{-5}$ & $1.30\cdot 10^{-4}$ & $1.60\cdot 10^{-4}$ \\
$x_0$ removal             & $9.78\cdot 10^{-5}$ & $1.28\cdot 10^{-4}$ & $1.55\cdot 10^{-4}$ \\
$x_0$ removal + full cond & $3.14\cdot 10^{-5}$  & $3.44\cdot 10^{-5}$ & $6.71\cdot 10^{-5}$ \\
$x_0$ removal + part cond & $3.80\cdot 10^{-5}$ & $4.28\cdot 10^{-5}$ & $6.72\cdot 10^{-5}$ \\
\bottomrule
\end{tabular}
\end{table}

\subsection{Solving different approximations of quadratic constraints}

In Table~\ref{tab:qp_qcqp:2}, there are the computation times for solving a baseline QCQP, respectively QPs approximating the same quadratic constraints with a different number of affine constraints.
More precisely, in this case there are two masses, a bounded force acts on the first mass while there is a softed quadratic constraint (per stage) on the total energy of the second mass.
There is no penalty on the states of the system, with the controller only aiming at reducing the second mass energy to within the (initially infeasible) constrained level.
All problems have $N=6$ horizon, $n_x=4$ states, $n_u=1$ controls, $n_{b_u}=1$ control bounds; additionally they have:
\begin{itemize}
\item QCQP\_$\infty$: 1 softed state quadratic constraint per stage
\item QP\_4: a softed square approximation (implemented using 2 two-sided softed state bounds) per stage 
\item QP\_6: a softed hexagon approximation (implemented using 3 two-sided softed general affine bounds) per stage 
\item QP\_8: a softed octagon approximation (implemented using 4 two-sided softed general affine bounds) per stage 
\end{itemize}
In order to focus on the effect of the different affine approximations, the problems are solved without (partial) condensing.
From the table, we can see that the QCQP with the exact quadratic formulation can be solved using the QCQP IPM solver in a time in between the cases of hexagon and octagon affine approximations solved using a QP IPM solver.
On the other hand, the square approximation implemented using the cheaper box constraints is more noticeably faster, but still only in the order of 20\%.

\begin{table}
\centering
\caption{Total solution times (in seconds, for fixed 9 IPM iterations plus $x_0$ removal) for the mass-spring QCQP and QPs with different affine approximations of the quadratic constraints}
\label{tab:qp_qcqp:2}
\begin{tabular}{@{}llll@{}}
\toprule
QCQP\_$\infty$ & QP\_4 & QP\_6 & QP\_8 \\
\midrule
$9.60\cdot 10^{-5}$ & $7.76\cdot 10^{-5}$ & $9.37\cdot 10^{-5}$ & $9.93\cdot 10^{-5}$ \\
\bottomrule
\end{tabular}
\end{table}

\subsection{Comparison with additional software}

\begin{figure}[tpb]
	\centering
	\includegraphics[width=\linewidth]{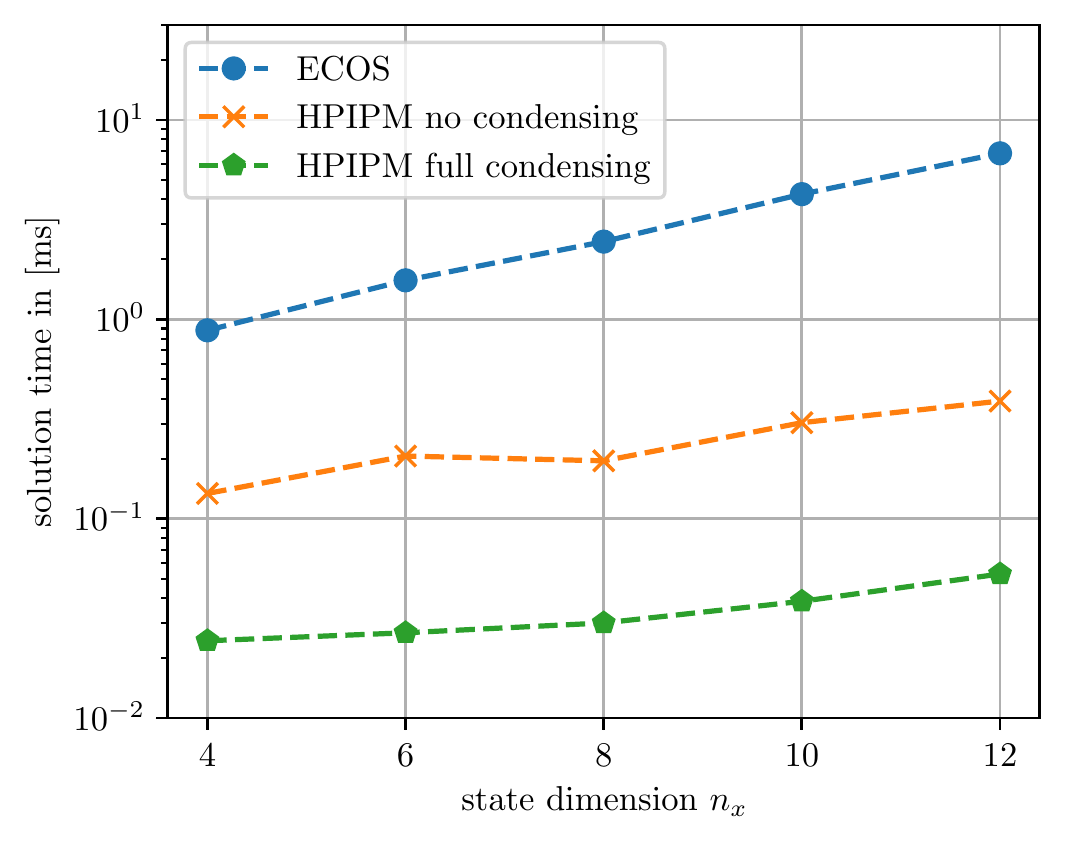}
	\caption{Computation time comparison ECOS and HPIPM.}
	\label{fig:eocs}
\end{figure}

We now compare the presented QCQP solver implementation within HPIPM with and without (full) condensing to ECOS, a sparse, open-source IPM for SOCP with a focus on embedded applications \cite{Domahidi2013a}.
The mass-spring OCP QCQP is solved to convergence for different number of masses, $n_u = 1$ controls and a horizon of $ N = 15 $.
The computation times are compared in Figure \ref{fig:eocs}.
The results show that using the presented efficient OCP QCQP framework can result in a computation time reduction of up to two orders of magnitude compared to an efficient state-of-the-art SOCP solver.
For $ n_x = 4 $, using HPIPM (without condensing) compared to ECOS results in a speed up of roughly factor 10.
For $ n_x = 12 $, using HPIPM with full condensing is approximately two orders of magnitude faster compared to ECOS.
One can observe that, compared to ECOS, the computation time of HPIPM grows slower when increasing the state dimension, since the flops count increase is partially offset by the increased performance of the underlying BLASFEO routines, which reach close to peak performance for matrices of sizes of multiple tens.
The performance of the BLASFEO routines also gives a typical stairs-like figure, following the granularity of its computational kernels (in this case of width 4).

%

\section{Conclusion \& Outlook}

We have shown a significant performance increase of the presented IPM tailored to QCQP compared to a state of the art SOCP solver for embedded optimization.
Moreover, the computational effort is of the same order of magnitude compared to solving a comparable QP.
When also the proposed efficient (partial) condensing algorithms for QCQP are taken into account, we are confident that the proposed QCQP framework in HPIPM can prove itself successful in many applications.
In that regard, future work includes a convenient interface into the higher-level MPC software framework \texttt{acados}, \cite{Verschueren2021}.
The efficient SQP methods therein can be extended to SQCQP, which would simplify a thorough comparison of the two when regarding real-world MPC tasks.


\bibliographystyle{abbrv}
\bibliography{syscop}

\end{document}